\documentclass[12pt]{article}

\usepackage{amsfonts,amsmath,amsthm,epsfig,verbatim,amssymb,amscd,eucal,bbm}

\usepackage[latin1]{inputenc}

\def\P{{\mathbb{P}}}
\def\E{{\mathbb{E}}}
\def\cG{{\cal{G}}}

\newtheorem{theorem}{Theorem}
\newtheorem{lemma}[theorem]{Lemma}

\begin{document}

%\begin{frontmatter}
%% "Title of the paper"
\title{Thresholds for virus spread on networks}
%\runtitle{Thresholds for virus spread on networks}

% indicate corresponding author with \corref{}
% \author{\fnms{John} \snm{Smith}\corref{}\ead[label=e1]{smith@foo.com}\thanksref{t1}}
% \thankstext{t1}{Thanks to somebody}
% \address{line 1\\ line 2\\ printead{e1}}
% \affiliation{Some University}

\author{ M. {\sc Draief} \thanks{Statistical Laboratory, Centre for
Mathematical Sciences, Wilberforce Road, Cambridge CB3 0WB UK
E-mail: {\tt M.Draief@statslab.cam.ac.uk}},~A. {\sc Ganesh}
\thanks{Microsoft Research, 7 J.J. Thomson Avenue, Cambridge CB3 0FB E-mail: {\tt ajg@micorsoft.com} } ~ and ~ L. {\sc Massoulié}
\thanks{Microsoft Research, 7 J.J. Thomson Avenue, Cambridge CB3 0FB E-mail: {\tt lmassoul@micorsoft.com} }}

\date{}

\maketitle

\begin{abstract}
We study how the spread of computer viruses, worms, and other
self-replicating malware is affected by the logical topology of the
network over which they propagate. We consider a model in which each
host can be in one of 3 possible states - susceptible, infected or
removed (cured, and no longer susceptible to infection). We
characterise how the size of the population that eventually becomes
infected depends on the network topology. Specifically, we show that
if the ratio of cure to infection rates is larger than the spectral
radius of the graph, and the initial infected population is small,
then the final infected population is also small in a sense that can
be made precise. Conversely, if this ratio is smaller than the
spectral radius, then we show in some graph models of practical
interest (including power law random graphs) that the final infected
population is large. These results yield insights into what the
critical parameters are in determining virus spread in networks.
\end{abstract}
%
%\begin{keyword}[class=AMS]
%\kwd[Primary ]{05C80} \kwd{90B15} \kwd[; secondary ]{68R99}
%\end{keyword}

%\begin{keyword}
%\kwd{Reed-Frost epidemic, random graphs, epidemic threshold,
%spectral radius, giant component} \kwd{}
%\end{keyword}
%%
%\end{frontmatter}

\section{Introduction}
Computer viruses and worms are self-replicating pieces of code
that propagate in a network. The essential difference between them
is that a virus typically needs some form of human intervention,
such as opening an email attachment or executing some software, to
cause them to be replicated, whereas worms do not require such
intervention. They use a number of different methods to identify
new targets for infection; for example, many worms scan randomly
generated IP addresses to locate vulnerable hosts, while email
viruses send copies of themselves to all addresses in the address
book of the victim. A survey of techniques for target location can
be found in \cite{taxonomy}.

The particular mechanism chosen by a worm or virus to propagate
itself defines a topology over which the infection can potentially
spread. What impact does the topology have on the speed of spread
of the epidemic, and moreover what are the key features of the
topology that determine how virulent the epidemic is? These are
the questions that we address in this paper.

In this paper, we use a susceptible-infected-removed (SIR) model
to describe the spread of the epidemic. Here, each susceptible
node can be infected by its infected neighbours at a rate
proportional to their number, and remains infected for a
deterministic or random time until it is removed. While it is
infected, it has the potential to infect its neighbours. Removal
can correspond to either (i) patching the computer represented by
the node, or (ii) its disconnection from the network by some
quarantining mechanism, or (iii) the exhaustion of its infectious
period either by a time-out mechanism or because it has tried all
its neighbours. Once a node is removed, it cannot become
susceptible or infected again. Our model ignores the possibility
that susceptible nodes can also be removed, e.g., because they
have received a patch or virus signature conferring immunity. This
is justified if the timescale for patching of susceptible hosts is
much larger (happens much more slowly) than that of epidemic
spread.

%We will primarily be interested in the latter mechanism,
%operating on the timescale of epidemic spread.
In the context of worms, there has recently been considerable
interest in automatic mechanisms for detecting whether hosts are
infected, and throttling or quarantining them; see, e.g.,
\cite{williamson}. There has also been work on automatic generation
of self-certifying alerts \cite{costa} which are equivalent to
patches. Thus, it is possible to view removal as happening on the
same time scale as infection. In the case of viruses, it takes
longer to generate virus signatures and update antivirus software,
but their spread is also slower. Hence, again, the model is not
unrealistic.

There is a substantial literature on the SIR model in epidemiology,
starting with the work of Kermack and McKendrick \cite{KM}. A
commonly used approach in early work was to approximate a stochastic
model by a deterministic one in a large population (law of large
numbers) limit. More recent work has considered stochastic aspects,
such as obtaining Poisson or normal limiting distributions for the
number of survivors; see, for example, \cite{BU,LU}. A key concept
in these studies is the basic reproductive number $R_0$, which
denotes the expected number of secondary infectives caused by a
single primary infective. If $R_0>1$, the infection spreads to some
sizeable fraction of the entire population; if $R_0<1$, then the
fraction eventually infected is close to zero. The concept of basic
reproductive number is easy to define with uniform mixing (i.e.,
when any infective can infect any susceptible equally easily) but it
is not clear how to apply it to general networks, where this number
could be different for every node. One approach is to consider
networks with special structure, where either nodes or links belong
to one of a small number of types. This is the approach taken, for
example, by Ball et al. \cite{ball}, who consider two-level models
where network links can belong to one of two types - (i) local,
e.g., within a household or (ii) global, between households.

In this paper, we obtain conditions for the number eventually
infected to be small, in arbitrary networks. Conversely, we obtain
conditions for the number of infected nodes to be large in some
specific network models of practical interest, including
Erd\H{o}s-R\'enyi and power law random graphs. The rest of the paper
is structured as follows. We introduce the epidemic spreading model
in Section \ref{sec:model}. Sufficient conditions for small epidemic
size (where the size is defined as the number that ever become
infected) are obtained in Section \ref{sec:model}. Applications of
these results to the star, clique, Erd\H{o}s-R\'enyi graph, and
power law graph are found in Section \ref{sec:applications}.
Furthermore in this Section we took advantage of results on the
giant component for various families of graphs to give a lower bound
to the number of nodes ultimately removed. Section \ref{sec:summary}
summarizes the paper and describes further directions to pursue.

\section{Model}\label{sec:model}
We consider a closed population of $n$ individuals, connected by a
neighbourhood structure which is represented by an undirected,
labelled graph $\mathcal{G}=(V,E)$ with node set $V$ and edge set
$E$. Each node can be in one of three possibly states, susceptible
(S), infective (I) or removed (R). The initial set of infectives
at time $0$ is assumed to be non-empty, and all other nodes are
assumed to be susceptible at time $0$. The evolution of the
epidemic is described by the following discrete-time model. Let
$X_v(t)$ denote the indicator that node $v$ is infected at the
beginning of time slot $t$ and $Y_v(t)$ the indicator that it is
removed. Each node that is infected at the beginnning of a time
slot attempts to infect each of its neighbours; each infection
attempt is successful with probability $\beta$ independent of
other infection attempts. Each infected node is removed at the end
of the time slot. Thus, the probability that a susceptible node
$u$ becomes infected at the end of time slot $t$ is given by
$1-\prod_{v\sim u}( 1-\beta X_v(t) )$, where we write $v\sim u$ to
mean that $(u,v)\in E$. Note that the evolution stops when there
are no more infectives in the population. At this time, we want to
know how many nodes are removed.

The above model is known as the Reed-Frost model. It corresponds
to a deterministic infectious period which is the same at every
node. It is one of the earliest stochastic SIR models to be
studied in depth, because of its analytical tractability. Note
that the evolution can be described by a Markov chain in this
case. Another commonly used model assumes that infectious periods
are iid and exponentially distributed, so that the system evolves
as a continuous time Markov process. General infectious periods
give rise to non-Markovian systems. These are outside the scope of
this work.

The object of interest is the number of nodes that eventually
become infected (and removed) compared to the number initially
infected. As noted earlier, in mean field models of SIR epidemics,
the number of nodes removed exhibits a sharp threshold; as $\beta$
is increased, it suddenly jumps from a constant (which doesn't
depend on $n$) to a non-zero fraction of $n$, the number of nodes
in the system. We wish to ask if a similar threshold is exhibited
on general graphs and, if so, how the critical value of $\beta$ is
related to properties of the graph.

We now state general conditions for the number of nodes removed to
be small. Let $A$ denote the adjacency matrix of the undirected
graph $G$, i.e., $a_{ij}=1$ if $(i,j)\in E$ and $a_{ij}=0$
otherwise. Since $A$ is a symmetric, non-negative matrix, all its
eigenvalues are real, the eigenvalue with the largest absolute
value is positive and its associated eigenvector has non-negative
entries (by the Perron-Frobenius theorem). If the graph is
connected, as we shall assume, then this eigenvalue has
multiplicity one, and the corresponding eigenvector is the only
one with all entries being non-negative.

\begin{theorem} \label{thm:upperbound_general}
Suppose $\beta \lambda_1<1$. Then, the total number of nodes
removed, $|Y(\infty)|$, satisfies
$$
\E [|Y(\infty)|] \leq \frac{1}{1-\beta \lambda_1} \sqrt{n |X(0)|},
$$
where $|X(0)|$ is the number of initial infectives. Morevoer, if
the graph $G$ is regular (i.e., each node has the same number of
neighbours), then
$$
\E[ |Y(\infty)| ] \le \frac{1}{1-\beta \lambda_1} |X(0)|.
$$
\end{theorem}

\begin{proof} In order for an arbitrary node $v$ to be infected at
the start of time slot $t$, there must be a chain of distinct
nodes $u_0, u_1,\ldots,u_t=v$ along which the infection passes
from some initial infective $u_0$ to $v$. Thus, by the union
bound,
$$
\P(X_v(t)=1) \le \sum_{u_0,\ldots,u_{t-1}} \beta^t X_{u_0}(0),
$$
where the sum is taken over nodes $u_0,\ldots,u_{t-1}$ such that
$(u_{i-1},u_i)\in E$ for all $i=1,\ldots,t$, where we take
$u_t=v$. Note that we have not imposed the requirement that the
$u_i$ be distinct as we are only seeking an upper bound.
Consequently, the probability that node $v$ ever gets infected
(and hence that $Y_v(\infty)=1$) is bounded above by
$$
\P(Y_v(\infty)=1) \le \sum_{t=0}^{\infty} \sum_{u\in V} (\beta
A)^t_{uv} X_u(0),
$$
since the $uv^{\rm th}$ entry of the matrix $A^t$ is simply the
number of paths of length $t$ between nodes $u$ and $v$. It is
immediate from the above that
$$
\E [|Y(\infty)|] = \sum_{v\in V} \P (Y_v(\infty)=1) \le
\sum_{t=0}^{\infty} {\bf 1}^T (\beta A)^t X(0),
$$
where ${\bf 1}$ denotes the vector of ones. Now, if $\beta
\lambda_1<1$, then we can rewrite the above as
\begin{eqnarray}
\E [|Y(\infty)|] &\le& {\bf 1}^T (I-\beta A)^{-1} X(0) \nonumber \\
&\le& \| {\bf 1} \| \ \| (I-\beta A)^{-1} \| \ \| X(0) \|,
\label{ubd1}
\end{eqnarray}
where $\| \cdot \|$ denotes the Euclidean norm in the case of a
vector, and the matrix or operator norm in the case of a matrix.
Now the operator norm of a symmetric matrix is its spectral
radius, the largest of its eigenvalues in absolute value. Hence
$\| (I-\beta A)^{-1} \| = (1-\beta \lambda_1)^{-1}$. Moreover, $\|
X(0) \| = \sqrt{\sum_{v\in V} X_v^2(0)} = \sqrt{ |X(0)| }$.
Likewise, $\| {\bf 1} \| =\sqrt{n}$. Substituting these in
(\ref{ubd1}) yields
$$
\E [|Y(\infty)|] \le \frac{1}{1-\beta \lambda_1} \sqrt{n |X(0)|},
$$
which is the first claim of the theorem.

Next, note that by using the spectral decomposition
$$
(I-\beta A)^{-1} = \sum_{i=1}^n \frac{1}{1-\beta \lambda_i} x_i
x_i^T,
$$
where $x_i$ denotes the eigenvector corresponding to the
eigenvalue $\lambda_i$ of $A$, and $x_i^T$ its transpose, we can
rewrite (\ref{ubd1}) as
\begin{equation} \label{ubd2}
\E [|Y(\infty)|] \le \sum_{i=1}^n \frac{1}{1-\beta \lambda_i} {\bf
1}^T x_i x_i^T X(0).
\end{equation}
Now, if $G$ is a regular graph and each node has degree $d$ (i.e.,
has exactly $d$ neighbours), then each row sum of its adjacency
matrix $A$ is equal to $d$. Hence, it is clear that the positive
vector $\frac{1}{\sqrt n}{\bf 1}$ is an eigenvector of $A$
corresponding to the eigenvalue $d$. By the Perron-Frobenius
theorem, this is therefore the largest eigenvalue. Hence,
$\lambda_1 = d$, $x_1 = \frac{1}{\sqrt n}{\bf 1}$, and all other
eigenvectors $x_2,\ldots,x_n$ are orthogonal to ${\bf 1}$. Hence,
by (\ref{ubd2}),
\begin{eqnarray*}
\E [|Y(\infty)|] &\le& \frac{1}{1-\beta \lambda_1} {\bf 1}^T x_1 x_1^T X(0) \\
&=& \frac{1}{n(1-\beta \lambda_1)} {\bf 1}^T {\bf 1} {\bf 1}^T X(0)
= \frac{1}{1-\beta \lambda_1} |X(0)|.
\end{eqnarray*}
This is the second claim of the theorem.
\end{proof}

Actually, there is an easier proof. Let $\nu(i) = \P(\mbox{node
$i$ is ever infected})$. Then, $\nu(i)=1$ if $i\in I$, where $I$
denotes the set of initial infectives, and otherwise $\nu(i) \le
\sum_{j\sim i} \beta \nu(j)$, where we write $j\sim i$ to mean
that $(i,j)$ is an edge. Thus,
\begin{equation} \label{infect_prob_bd}
(I-\beta A) \nu \le {\bf 1}_I,
\end{equation}
where ${\bf 1}_I$ denotes the vector with components 1 for $i\in
I$ and 0 for $i\notin I$, and the inequality holds in the usual
partial order, namely componentwise. Now, if $\beta
\lambda_1(A)<1$, then we have the power series expansion
$$
(I-\beta A)^{-1} = \sum_{k=0}^\infty \beta^k A^k,
$$
from which it is immediate that $(I-\beta A)^{-1}$ is a
non-negative matrix. Therefore, we can multiply both sides of the
inequality in (\ref{infect_prob_bd}) by $(I-\beta A)^{-1}$ to
obtain
$$
\nu = \E[Y(\infty)] \le (I-\beta A)^{-1} X(0),
$$
and so
$$
\E[|Y(\infty)|] \le {\bf 1} (I-\beta A)^{-1} X(0).
$$
This is the same as (\ref{ubd1}), and the proof carries on the
same way from there.

\noindent{\bf Remarks} The upper bound in the first claim of the
theorem is close to the best possible in general, as the example
of the star-shaped network in Section \ref{sec:star} demonstrates.

The theorem says that, if $\beta \lambda_1 < 1$, then starting
from a `small' population of initial infectives, the final size of
the epidemic is small. For example, if $|X(0)|=1$, then the final
size of the epidemic is bounded by a constant in the case of
regular graphs, and by a multiple of $\sqrt{n}$ in general. Thus,
the fraction of nodes infected goes to zero as $n$ tends to
infinity.

Note that the proof of the theorem above doesn't require us to
assume that the epidemic be of Reed-Frost type. It works for
general infectious periods $J$ since we are only using
expectations throughout, which don't require independence
assumptions. Therefore, following the steps of the above proof and
replacing $\beta$ by the probability that a node gets infected bu
an infected neighbour.
%Let $J_u$ be the length of the infectious
%period of node $u$ and let $E_v$ the time needed to infect its
%neighbour $v$, $(E_v)_{v\in V}$ form an i.i.d sequence of
%exponentially distributed random variables with parameter
%$\lambda$.
%$$P(X_u(t)=1)\led \sum_{u\sim v}X_u(t-1)P(J_u<.

In turn, if node $u$ is infected, it will infect $j$, if they are
connected and if the time it takes to contact this node given by
an exponential random variable with parameter $\lambda$ is less
than $J$.

\begin{theorem} \label{thm:upperbound_general}
Suppose that $J$ is such that $\E[e^{-\lambda J}]<\infty$ and let
$$p_J=1-\E[e^{-\lambda J}]\:.$$ If $p_J \lambda_1<1$ then the total
number of nodes removed, $|Y(\infty)|$, satisfies
$$
\E [|Y(\infty)|] \le \frac{1}{1-p_J \lambda_1} \sqrt{n |X(0)|},
$$
where $|X(0)|$ is the number of initial infectives. Morevover, if
the graph $G$ is regular (i.e., each node has the same number of
neighbours), then
$$
\E[ |Y(\infty)| ] \le \frac{1}{1-p_J \lambda_1} |X(0)|.
$$
\end{theorem}

The converse is not true in general. Consider the ring on $n$
nodes. As a regular graph with degree 2, its adjacency matrix has
maximum eigenvalue 2. Nevertheless, for any $\beta<1$, the size of
the epidemic starting from a single initial infective is bounded
by the sum of the sizes of two independent branching processes,
where each branching process has Bernoulli offspring distribution
with parameter $\beta$. (The branching processes decide if the
epidemic will spread one node left or right from its current
position before dying out). Since each branching process is
subcritical, its final size is finite almost surely and in
expectation. Thus, the expected size of the epidemic is a constant
that does not depend on $n$. In other words, the epidemic is small
even if $\beta \lambda_1 = 2\beta > 1$.

In particular, the SIR epidemic on the ring does not exhibit a
sharp threshold on the open interval $(0,1)$. On this interval,
the final size of the epidemic is a smooth function of the
infectiousness parameter $\beta$, even in the limit as $n$ tends
to infinity. It is shown in Section \ref{sec:star} that a similar
result holds for star-shaped networks as well; in fact, there is
no threshold even on the closed interval $[0,1]$ in this case.

However, while there isn't always a threshold, it turns out that
there is one in many networks of practical interest: there is a
lower bound on $\beta$, above which the epidemic infects a
positive fraction of the population on average. We now illustrate
this through several examples.

\section{Examples} \label{sec:applications}

\subsection{Star-shaped networks} \label{sec:star}

The star-shaped network is of interest because it illustrates that
the bound in Theorem \ref{thm:upperbound_general} is close to the
best possible for general networks. It also exhibits a smooth
dependence of the final size of the epidemic on the infectiousness
parameter $\beta$, thereby demonstrating that threshold behaviour
doesn't always occur. Finally, understanding the star is important
to understanding certain power-law networks.

Consider the star network, consisting of a hub and $n-1$ leaves,
each of which is attached only to the hub. Its adjacency matrix
$A$ has ones along the first row and column, except for the
$(1,1)$ element, which is zero; all other elements are zero. In
other words, $A = {\bf 1}{\bf 1}^T - e_1 e_1^T$, where $e_1^T = (1
0 \ldots 0)$. Thus $A$ is a rank-two matrix and can have only two
non-zero eigenvalues. It is readily verified that $(\sqrt{n-1} 1
\ldots 1)^T$ and $(-\sqrt{n-1} 1 \ldots 1)^T$ are eigenvectors
corresponding to the eigenvalues $\sqrt{n-1}$ and $-\sqrt{n-1}$
respectively, and so the spectral radius of $A$ is $\sqrt{n-1}$.

Now suppose $\beta \sqrt{n-1}=c<1$. Consider the initial condition
where only the hub is infected, so that $|X(0)|=1$. The number of
leaves infected before the hub is cured is binomial with
parameters $n-1$ and $\beta$. No other leaves can be infected
subsequently. Hence,
$$
\E [|Y(\infty)|] = 1+\beta(n-1) = 1+c\sqrt{n-1},
$$
which is comparable to the upper bound, $\sqrt{n-1}/(1-c)$, given
by Theorem \ref{thm:upperbound_general}. We also observe in this
case that $\E[|Y(\infty)|]$ is a smooth (almost linear) function
of $\beta$ and does not exhibit any threshold behaviour.

Suppose next that the hub is initially uninfected but $k$ leaves
are infected. The hub becomes infected in the next time step with
probability $1-(1-\beta)^k$. It subsequently infects a number of
leaves which is binomial with parameters $n-1-k$ and $\beta$. The
epidemic dies out at $t=3$. So, in this case,
\begin{eqnarray*}
\E [|Y(\infty)|] &=& k + [1 - (1-\beta)^k][1+\beta(n-1-k)] \\
&\le& k + \beta k[1+\beta(n-1-k)] \le |X(0)| (1+2c^2).
\end{eqnarray*}
Thus, when the hub is initially uninfected, the expected final
size of the epidemic is only a constant multiple of the initial
number of infectives. This illustrates that the initial condition
can have a big impact in general.

\subsection{Complete graph}
A complete graph is one which an edge is present between every
pair of nodes. Much of the early work on SIR epidemics was based
on mean field models.  These are rigorously justifiable only in
the case of complete graphs, and motivates our interest in them.
We shall recover the classical result that the epidemic has a
threshold at $R_0=1$, where the basic reproduction number $R_0 =
\beta(n-1)$ is defined as the mean number of secondary infections
caused by a single primary infective, when the entire population
is susceptible. From the perspective of networking applications,
the BGP routers belonging to the top level autonomous systems of
the Internet form a completely connected component. In addition,
large ISPs often organize their internal BGP (iBGP) routers into a
set of \emph{route reflectors} that are completely connected.

The complete graph is a regular graph with common node degree
$n-1$. Therefore, its spectral radius is $\lambda_1 = n-1$, and we
have by Theorem \ref{thm:upperbound_general} that, if
$\beta(n-1)<1$, then the final size of the epidemic is bounded by
$1/(1-\beta(n-1))$ times the initial number of infectives. We now
establish a converse.

Suppose $\beta(n-1) = c>1$ is held constant. (We don't need to
assume this, but the results will need to be restated in terms of
the limits superior and inferior of the sequence $c_n$; it should
be clear to the reader how to do so based on the discussion
below.) Let $|X_0|=1$ and let $u$ be the initial infected node.
Consider the random subgraph of the complete graph obtained by
retaining each edge with probability $\beta$, independent of all
other edges, and let $C_u$ denote the connected component
containing $u$ in this random graph (possibly just the singleton
$\{ u \}$). It is clear that $C_u$ can be interpreted as the set
of infected nodes in the epidemic; each neighbour of $u$ is
infected with probability $\beta$, and is hence a neighbour of $u$
in the random graph described above, and so on iteratively. Thus,
the number of infected nodes in the epidemic has the same
probability law as the size of the component $C_u$.

The above random graph model was introduced by Erd\H{o}s and
R\'enyi \cite{ER60}; we denote it by $\cG(n,\beta)$, where $n$
denotes the number of nodes, and $\beta$ the probability that the
edge between each pair of nodes is present. It is also called a
Bernoulli random graph because the indicators of edges are iid
Bernoulli random variables.

We now use the following fact, which was established by Erd\H{o}s
and R\'enyi \cite{ER60}; see \cite[Theorem 5.4]{JaLuRu00}, for
instance, for a more recent reference. Here, we assume that
$c=\beta(n-1)$ is held constant while $n\to \infty$, and that
$c>1$.

\begin{theorem} \label{thm:giant}
Let $\gamma$ be the unique positive solution of $\gamma +
e^{-\gamma c}=1$. Then, as $n\to \infty$, the size of the largest
connected component in the random graph $\cG(n,\beta)$ is
$(1+o(1))\gamma n$, with probability going to $1$ as $n$ tends to
infinity.
\end{theorem}

The uniqueness of $\gamma$ follows from the convexity of the
function $f(x)=x+e^{-cx}$ and the fact that $f(0)=1$, while its
existence follows from the continuity of $f$ and the fact that
$f'(0) = 1-c < 0$, but that $f(x) \to \infty$ as $x\to \infty$.

We now estimate the size of $C_u$, the connected component
containing the initial infective. If $u$ belongs to the `giant
component', then $|C_u| \equiv \gamma n$. Since a fraction
$\gamma$ of nodes belong to the giant component, the probability
that node $u$ does so is $\gamma$. Hence, $\E [|C_u|] =
(1+o(1))\gamma^2 n$. We have thus shown the following:

\begin{lemma} \label{lem:complete}
Let $G=(V,E)$ be the complete graph on $n$ nodes, and let $\beta =
\frac{c}{n-1}$ for an arbitrary constant $c>1$. Then, the final
size of the epidemic satisfies
$$
\E [|Y(\infty)|] \ge (1+o(1)) \gamma^2 n,
$$
for any $|X(0)| \ge 1$, where $\gamma>0$ solves $\gamma+e^{-\gamma
c}=1$.
\end{lemma}

There is thus a threshold at $c=1$ for the final size of the
epidemic; starting with a constant number of initial infectives,
the final size is a constant independent of $n$ if $c<1$, and a
fraction of $n$ if $c>1$.

\subsection{Erd\H{o}s-R\'enyi random graphs}
The Erd\H{o}s-R\'enyi graph $\cG(n,p)$ with parameters $n$ and $p$
is defined as a random graph on $n$ nodes, where the edge between
each pair of nodes is present with probability $p$, independent of
all other edges. If $p=1$, then this is the complete graph.

The spreading behavior of an epidemic on an
Erd\H{o}s-R\'enyi graph is of interest for a number of reasons. First,
it is a graph that has received considerable attention in the past
\cite{bollobas}.  Second, it is an important component of the class
of power law random graphs that model the Internet AS graph. Thus if we
are to understand the robustness of the Internet AS-level graph, we
need to characterize the robustness of the Erd\H{o}s-R\'enyi graph.

We shall consider a sequence of such graphs indexed by $n$. Denote
by $d$ the corresponding average degree, i.e. $d=(n-1)p$. Note
that $p$ and $d$ depend on $n$, but this is suppressed in the
notation. We say that a property holds with high probability if
its probability goes to 1 as $n\to \infty$. Define $c_n =\beta d =
(n-1)\beta p$; we have suppressed the dependence of $\beta$ and
$p$ on $n$ in the notation, but make it explicit in the case of
$c$. Consider an SIR epidemic on such a graph starting with one
node initially infected. We have the following:

\begin{lemma} \label{lem:er}
If $\limsup_{n\to \infty} c_n \le c<1$, then for all $n$
sufficiently large, $\E[ |Y(\infty)| ]$ is bounded by a constant
that does not depend on $n$. On the other hand, if $\liminf_{n\to
\infty} c_n \ge c > 1$, then $\E [|Y(\infty)| ] \ge
(1+o(1))\gamma^2 n$ where $\gamma>0$ solves $\gamma+e^{-\gamma c}
= 1$.
\end{lemma}

\begin{proof} Suppose first that $\liminf_{n\to \infty} c_n \ge
c>1$. As in the case of the complete graph, we identify the
infected individuals in the epidemic with the connected component
containing the initial infective $u$ in an Erd\H{o}s-R\'enyi
random graph with parameters $n$ and $\beta p$. (If edge $(u,v)$
is present in the original Erd\H{o}s-R\'enyi graph, which happens
with probability $p$, then $u$ succeeds in infecting $v$ with
probability $p$. This yields the new graph with edge probability
$\beta p$; the independence of the edges is obvious.) Thus, the
second claim of the lemma follows in the same way as Lemma
\ref{lem:complete}.

The first claim is stronger than what the upper bound of Theorem
\ref{thm:upperbound_general} yields. Note that, by the
Perron-Frobenius theorem, the spectral radius $\lambda_1(A)$ of
the adjacency matrix lies between the smallest and largest node
degree.  For the random graph $\cG(n,p)$, the node degrees are
binomial random variables with parameters $n-1$ and $p$. If the
average node degree $d=(n-1)p$ satisfes $d \gg \log(n)$, i.e.,
$\log(n)/d \to 0$ as $n\to \infty$, then it can be shown using
Chernoff's bound that both the minimal and maximal node degree are
$(1+o(1))d$ with high probability; hence, so is the spectral
radius. In this case, Theorem \ref{thm:upperbound_general} yields
that, if $\beta \lambda_1(A) \sim (n-1)\beta p \le c<1$, then the
expected final size of the epidemic is bounded by a constant times
$\sqrt{n}$. To show that it is in fact bounded by a constant, and
that this holds even without the assumption that $d \gg \log(n)$,
we use a branching process construction.

Rather than fixing the random graph $\cG(n,p)$ in advance, we use
the principle of deferred decisions to generate it dynamically as
the epidemic progresses. Thus, starting with the initial infective
$u$, we put down all edges from it to other nodes. Then, we decide
whether $u$ succeeds in infecting its neighbours along each of
those edges. For each neighbour $v$ so infected, we repeat the
process. Thus, the number of nodes infected by $u$ is binomial
with parameters $n-1$ and $\beta p$; the number of nodes infected
by each subsequent infective is stochastically dominated by such a
binomial random variable. Thus, the size of the epidemic is
bounded above by the size of a branching process whose offspring
distribution is binomial, $B(n-1,\beta p)$. The branching process
is subcritical by the assumption that $(n-1)\beta p \le c < 1$,
and so it becomes extinct with probability $1$, i.e., its final
population size is finite almost surely, and in expectation. It
can be shown directly, using generating functions, that it is
bounded uniformly in $n$. Alternatively, note that if $(n-1)\beta
p = c$ for all $n$, then the binomial offspring distributions
converge in distribution to a Poisson with parameter $c$ as $n\to
\infty$; the population sizes of the corresponding branching
processes also converge, both in distribution and in expectation.
Since $c<1$, the branching process with Poisson($c$) offspring
distribution is subcritical, and so it has a finite mean
population size. This completes the proof of the lemma.
\end{proof}

%\subsection{Hypercubes}
%The hypercube is of interest because of the widespread and growing
%interest in distributed hash tables and applications, such as file
%sharing \cite{pastry}, being built on top of them.  Already worms
%and viruses have appeared in some applications, \cite{kazaa}. As
%many DHT structures are hypercubic in nature, it is important to
%understand the spreading behavior of such worms on a hypercube. Here
%we represent a hypercube as a graph $G$ with vertex set
%$\{0,1\}^{\ell}$ for some $\ell\in\setN$, and where the edge $(v,w)$
%is present if and only if the Hamming distance $d_H(v,w)$ equals 1.
%As a hypercube is a regular graph, its spectral radius is
%$\lambda_1(G) = \log_2 n = \ell$. Hence, we have by Theorem
%\ref{thm:upperbound_general} that, if $\beta \ell < 1$, then the
%final size of the epidemic is bounded by a constant multiple of
%the number of initial infectives; the constant is $1/(1-\beta \ell)$,
%and does not depend on $n$.

\subsection{Power law random graphs}\label{sec:powerlaw}
There has been considerable interest in power law graphs since it
was first noticed that the Internet AS-level graph exhibits a
power law degree distribution, \cite{FFF}.  Briefly a power law
graph is one where the number of nodes with degree $k$ is
proportional to $k^{-\gamma}$ for some $\gamma > 1$.  For the mean
degree to be finite, we need $\gamma > 2$ and this is the range we
shall consider. The Internet AS-level graph is characterized by
$\gamma \approx 2.1$.

There have been several different models proposed for graphs with
power law degree distributions; see, for example, \cite{BA,BR}. In
this paper, we consider the following model of random graphs on
$n$ vertices, introduced in \cite{ChLu03}. Let ${\bf w} =
(w_{1},w_2,\ldots,w_{n})$ be a sequence of positive weights
assigned to the nodes of the graph; we assume without loss of
generality that $w_1 \ge$ $w_2 \ge \cdots $ $\ge w_n$. The edge
between the pair of vertices $(i,j)$ is present with probability
$$
p_{ij} = \frac{w_i w_j}{\sum_{k=1}^nw_k},
$$
independent of all other edges; we assume that $w_1^2 \leq
\sum_{k=1}^n w_k$. The resulting random graph is denoted $G({\bf
w})$. For example, taking $w_i=np$ for all $i$ yields the
Erd\H{o}s-R\'enyi model with parameters $n$ and $p$.

It is easy to see that $w_i$ is the expected value of the degree
of node $i$; hence, this model is referred to as the expected
degree model. \footnote{Reed and Molloy
\cite{molloyreed1,molloyreed2} have studied a model where one
conditions on actual rather than expected degrees. The expected
degree model has the advantage that edges are independent, which
makes it much easier to analyse.} We do not assume that the $w_i$
are integer-valued. Note that the resulting graph may have
self-loops but it does not have multiple edges. The self-loops do
not affect the spread of the epidemic and are not important to our
analysis.
%The presence of self-loops does not affect our conclusions;
%events that hold whp in this model will also hold whp in the
%same model conditioned not to have self-loops.

Let $d$ denote the average and $m$ the maximum expected degree.
(Thus $m=w_1$ but it is convenient to distinguish it in the
notation as the model is parametrised by $d$, $m$ and the exponent
of the power law degree distribution.) Chung and Lu \cite{ChLu03}
propose the following explicit power law model for the expected
degree sequence:
\begin{equation} \label{plrg_weights}
w_i = c(i_0+i)^{-\frac{1}{\gamma - 1}}, \; 1\leq i\leq n,
\end{equation}
where
\begin{equation} \label{plrg_constants}
c = \frac{\gamma - 2}{\gamma - 1}dn^{\frac{1}{\gamma -1}}, \; i_0
= n \Bigl( \frac{ d(\gamma -2)}{ m(\gamma - 1) } \Bigr)^{\gamma
-1}.
\end{equation}
The number of nodes with weight bigger than $k$ (equal to the
largest $i$ such that $w_i \ge k$) scales like $k^{1-\gamma}$.
Thus, $\gamma$ is the exponent of the power law distribution of
expected node degrees. The weights $w_i$ are the order statistics
of this distribution. The distribution is shifted by $i_0$ and
scaled by $c$, where these constants are chosen so as to achieve
the specified average $d$ and maximum $m$ for the expected
degrees.

The eigenvalues of the adjacency matrix for this random graph
model have been studied by Chung, Lu and Vu. They show
\cite[Theorem 4]{CLV04} that, with high probability, the spectral
radius of the graph is
$$\rho(A) = \left\{%
\begin{array}{ll}
(1+o(1))\sqrt{m}\:,  & \gamma > 2.5, \\
 (1+o(1))\frac{d(\gamma - 2)^2}{(\gamma - 1)(3 - \gamma)}
 \Bigl( \frac{(\gamma-1)m}{(\gamma - 2)d} \Bigr)^{3-\gamma}\:,
 & 2< \gamma < 2.5\:.
\end{array}%
\right.$$

By Theorem \ref{thm:upperbound_general}, if $\beta \rho(A)<1$,
then the size of the epidemic is bounded by $\sqrt{n}$ times the
size of the initial infective population.

%Is there a converse, i.e., is there a large epidemic (with
%positive probability) whenever $\beta \rho(A)>1$?

%We are not able to answer this question

%{\bf The rest of this section is poorly written. It contains many
%of the relevant calculations but we need to think about how to
%present it.}
%
%It turns out that the answer depends on the exponent $\gamma$
%of the power law degree distribution. The answer is yes if
%$\gamma<3$. If $\gamma>3$, then there is a large epidemic if
%the set of initial infectives contains a high-degree node but
%not otherwise; this is analogous to the star network studied earlier,
%where there is a large epidemic if the hub is initially infected, but
%otherwise the probability of a large epidemic is small. We make these
%statements precise below.

We now establish a partial converse. We show that the graph has a
core such that, if $\beta \rho(A)>1$ and one of the initial
infected nodes is in the core, then the expected size of the
epidemic is large. This is analogous to the situation in the star,
where there is a large epidemic when $\beta \rho(A)>1$ if the hub
is initially infected. For general power law graphs, we do not
know what happens when the initial infectives aren't in the core.

%(To put it another way, we cannot precisely characterise the
%core.)

%There have been several different models proposed for graphs with
%power law degree distributions; see, for example, \cite{BA,BR}.
%In this paper, we consider the following class of random graphs
%on $n$ vertices introduced in \cite{ChLu03}. It is defined by
%the sequence of expected degrees $(w_1,\dots,w_n)$ where vertex
%$i$ is assigned vertex weight $w_i$. The edges are chosen independently
%and randomly as follows: the probability that there is an edge
%between $i$ and $j$ is given by $p_{i,j}$ proportional to $w_iw_j$.
%In addition we assume that
%%$$\max_{i=1\dots n}w_i^2<w=\sum_{1}^nw_k\:,$$
%so that $p_{ij}\leq 1$, and the sequence $(w_1,\dots,w_n)$ is
%graphical\footnote{necessary and sufficient condition for a
%sequence to be realised by a graph}. We do not suppose the $w_i$'s
%to be integers. Note that the resulting graph may have self-loops
%but it does not have multiple edges. The self-loops do not affect
%the spread of the epidemic and are not important to our analysis.
%By way of example, the case $w_i=np$ for all $i$ corresponds to an
%Erd\H{o}s-R\'enyi random graph (modified to allow self-loops).

It is easy to see from the description above that the expected
degree of node $i$ is precisely $w_i$, and the expected average
degree of the graph is given by $d=\frac{1}{n}\sum_{i=1}^n w_i$.
We can now write $p_{ij}=\frac{w_iw_j}{nd}$.

It is straightforward to describe the evolution of a Reed-Frost
epidemic on the expected degree random graph model. Consider a
single initial infective, say node $i$. Node $j$ becomes infected
at time $1$ if edge $(i,j)$ is present in the random graph and if
$i$ infects $j$ in the first time slot; this has probability
$\beta p_{ij}$, and is independent of whether node $i$ infects
some other node $k$. Moreover, node $i$ cannot infect node $j$ in
any subsequent time step since it is removed at time $1$. Using
the principle of deferred decisions, we can construct a
realisation of the random graph as the epidemic spreads. It is
clear from this construction that the set of nodes that eventually
become infected can be identified with the connected components
containing the initial infectives in the random graph with weight
sequence $\beta {\bf w}$, namely $G(\beta {\bf w})$. Suppose there
is a single initial infective. The question of whether there is a
large epidemic is equivalent to that of whether the random graph
$G({\bf w})$ possesses a giant component, and whether the initial
infective belongs to this giant component. If there is more than
one initial infective, the final set of removed nodes is the union
of the connected components containing the initial infectives, in
the random graph $G({\bf w})$.

%The question of whether or not a large outbreak occurs is, as
%illustrated above, is closely related to the existence or not of a
%giant component.
A sufficient condition for the existence of a giant component is
derived in \cite[Theorem 3]{ChLu02}. The condition can be stated
in terms of the average expected degree $d$, as follows:
\begin{theorem} \label{thm:plrg_giant}
For a random graph $G(w)$ with expected degree sequence having
average expected degree $d>1+\delta>1$, there is a unique giant
component $C$ such that $\sum_{i\in C} w_i \ge (1-c_{\delta})
\sum_{i\in V} w_i$, where $c_{\delta} \in (0,1)$ is a constant
that depends only on $\delta$.
\end{theorem}
In words, the giant component contains a non-zero fraction of the
total weight of all nodes. Later, we will show that this implies
that it contains a non-zero fraction of the nodes.

We use this result to obtain estimates on the final size of an
epidemic on a power law random graph. Fix $k$ (as a function of
$n$) and consider the subgraph induced by the $k$ nodes with the
largest weight in the random graph $G(\beta {\bf w})$. The average
expected degree of this subgraph is easily seen to be
\begin{equation} \label{subgraph_degree}
d_k = \frac{ (\sum_{i=1}^k \beta w_i)^2 }{ k\sum_{i=1}^n \beta w_i
} = \frac{ \beta (\sum_{i=1}^k w_i)^2 }{ k nd}.
\end{equation}
If this is strictly larger than 1, then by Theorem
\ref{thm:plrg_giant} above, this subgraph has a giant component.
We now find conditions on $k$ such that $d_k > 1$.

We have from (\ref{plrg_weights}) and (\ref{plrg_constants}) that
$$
w_i = m \Bigl( 1+\frac{i}{i_0} \Bigr)^{-\frac{1}{\gamma-1}},
$$
and so,
\begin{eqnarray*}
\sum_{i=1}^k w_i \sim \int_0^k m \Bigl( 1+\frac{x}{i_0} \Bigr)^{-
\frac{1}{\gamma-1}} dx
&=& mi_0 \int_1^{1+(k/i_0)} y^{-\frac{1}{\gamma-1}} dy \\
&=& \frac{\gamma-1}{\gamma-2} mi_0 \Bigl[ \Bigl( \frac{k}{i_0}+1
\Bigr)^{\frac{\gamma-2}{\gamma-1}} - 1 \Bigr].
\end{eqnarray*}
Now, if $k \gg i_0$, it follows that
\begin{equation} \label{sumweights}
\sum_{i=1}^k w_i \sim \frac{\gamma-1}{\gamma-2} m
k^{\frac{\gamma-2}{\gamma-1}} i_0^{\frac{1}{\gamma-1}} = d
n^{\frac{1}{\gamma-1}} k^{\frac{\gamma-2}{\gamma-1}}.
\end{equation}
Substituting this in (\ref{subgraph_degree}) yields
\begin{equation} \label{subgraph_degree2}
d_k = \beta d n^{\frac{3-\gamma}{\gamma-1}}
k^{\frac{\gamma-3}{\gamma-1}}.
\end{equation}

We now distinguish two cases. Suppose first that $\gamma \ge 3$.
Then $d_k$ is a non-decreasing function of $k$, and its maximum
value, attained at $k=n$, is $\beta d$. This only yields the weak
result that there is a large epidemic if $\beta d>1$.
%On the other hand, simply by considering the star centred on the highest
%degree node, we note that if this node (or any of the $O(1)$ highest
%degree nodes) is initially infected, then the expected size of the epidemic
%is at least $\beta m$.

Suppose next that $2<\gamma<3$. Then $d_k$ is a decreasing
function of $k$. Fix $\delta>0$. Defining $N_{\delta}$ to be the
largest value of $k$ for which $d_k > 1+\delta$, we see that
\begin{equation} \label{coresize1}
N_{\delta} = \Bigl\lfloor \Bigl( \frac{\beta d}{1+\delta}
\Bigr)^{\frac{\gamma-1}{3-\gamma}} n \Bigr\rfloor + 1,
\end{equation}
where $\lfloor x \rfloor$ denotes the integer part of $x$. The
following result is now an easy consequence.

\begin{lemma} \label{threshold}
Let $\beta>0$ be arbitrarily small. Then the expected size of the
epidemic, starting from an arbitrary initial infective, is bounded
below by a constant multiple of $n$, where the constant may be
depend on $\beta$. Here the expectation is taken both with respect
to the random realisation of a graph, and the evolution of the
epidemic conditional on the graph.
\end{lemma}

\begin{proof}
Consider the $N_{\delta}$ nodes of largest weight, where
$N_{\delta}$ is given by (\ref{coresize1}). Since $\beta, d>0$ are
constants, $N_{\delta}$ is a constant multiple of $n$. By Theorem
\ref{thm:plrg_giant}, the random graph $G(\beta {\bf w})$
restricted to these nodes contains a giant component $C$ such that
\begin{equation} \label{giantweight}
\sum_{i\in C} w_i \ge (1-c_{\delta}) \sum_{i=1}^{N_{\delta}} w_i
\sim dn (1-c_{\delta}) \Bigl( \frac{\beta d}{1+\delta}
\Bigr)^{\frac{\gamma-2}{3-\gamma}},
\end{equation}
where we have used (\ref{sumweights}) and (\ref{coresize1}) to
obtain the last asymptotic equivalence. Recall that $c_{\delta}
\in (0,1)$ is a constant that depends on $\delta$. Since $\beta$,
$d$ and $\delta$ are constants, while $\sum_{i=1}^n w_i = nd$,
equation (\ref{giantweight}) tells us that the giant component $C$
contains a constant fraction of the total weight of the graph. We
now deduce that it must also contain a constant fraction of the
total number of nodes. Indeed, for a given weight, the size (in
number of nodes) would be minimised if $C$ contained the highest
weight nodes. Thus, we ask what is the smallest value of $k$ such
that $\sum_{i=1}^k w_i$ exceeds the weight of $C$. It follows from
(\ref{sumweights}) and (\ref{giantweight}) that we require
$$
dn^{\frac{1}{\gamma-1}} k^{\frac{\gamma-2}{\gamma-1}} \sim dn
(1-c_{\delta}) \Bigl( \frac{\beta d}{1+\delta}
\Bigr)^{\frac{\gamma-2}{3-\gamma}},
$$
and so,
\begin{equation} \label{giantsize}
k \sim n (1-c_{\delta})^{\frac{\gamma-1}{\gamma-2}} \Bigl(
\frac{\beta d}{1+\delta} \Bigr)^{\frac{\gamma-1}{3-\gamma}};
\end{equation}
in particular, $k$ is equivalent to a constant multiple of $n$.

Now, if any of the initially infected nodes belongs to the giant
component $C$, then $C$ is a subset of the set of nodes ever
infected; hence, the final size of the epidemic is proportional to
$n$. On the other hand, suppose none of the initial infectives
belongs to $C$. Let $i$ be an initially infected node. Now, the
probability that there is an edge between $i$ and $C$ in the
random graph $G(\beta {\bf w})$ is given by
$$
p(i,C) = \frac{w_i \sum_{j\in C} w_j}{nd} = (1+o(1))
(1-c_{\delta}) \Bigl( \frac{\beta d}{1+\delta}
\Bigr)^{\frac{\gamma-2}{3-\gamma}} w_i,
$$
which is a positive constant bounded away from zero. Conditional
on this edge being present, $C$ is a subset of the set of
eventually infected nodes. Thus, in this case too, the expected
final size of the epidemic is proportional to $n$. This completes
the proof of the lemma.
\end{proof}

Let us summarise our findings: If $\gamma>3$, then there is a
large epidemic if the set of initial infectives contains a
high-degree node but not otherwise; this is analogous to the star
network studied earlier, where there is a large epidemic if the
hub is initially infected, but otherwise the probability of a
large epidemic is small. In the next section, we introduce another
family of scale-free networks which will lead to more consistent
results.

\subsection{Inhomogeneous $W$-graphs}
We are interested in the following family of graphs $G(W)$ with
vertex set $\{1,\dots,n\}$ where $i$ and $j$ are connected by an
edge with probability $p_{i,j}$ independently from all other pairs
of nodes, and $p_{i,j}$ is defined as follows: let
$\{X_i,\:i=1,\dots,n\}$ a sequence of iid random variables
uniformly distributed on $[0,1]$ than $p_{ij}=W(X_i,X_j)$, where
$W:[0,1]^2\rightarrow[0,1]$ is measurable and symmetric. This
model has been introduced in \cite{LS04} and generalised in
\cite{BJR05} to a larger family of inhomogeneous networks. As a
matter of example, $W(x)=p$ for all $x$ gives the classical
Erdös-Rényi random graph. For more details on this family of
graphs we refer the reader to \cite{BJR05}. Our main interest in
what follows is to understand the emergence of the giant
component.

Let $T_W$ be the integral operator with kernel $nW$ defined by
$$T_Wf=\int_0^1 n W(x,y) f(y) dy\:,$$
and define
$$||T_W||=\sup\{||T_Wf||_2\: : \: f\geq 0,\: ||f||_2\leq 1\}\:.$$
As previously, we check that starting with a $G(W)$ graph, the
epidemic graph is described by $G(\beta W)$.

The next Theorem \cite[Corollary 3.2]{BJR05} characterises the
emergence of the giant component in this family of graphs.

\begin{theorem} Consider the graph
$G(W)$ then the threshold for the existence of the giant component
is given $1$. More precisely
\begin{itemize}
\item[$\bullet$] if $||T_W||\leq 1$, then the size of the largest
component is negligible with respect to $n$ almost surely,
\item[$\bullet$]  if $||T_W||>1$ and $W$ is irreducible\footnote{It cannot be split into two disjoint components}, then the largest
component consists of a non zero fraction of $n$.
\end{itemize}
\end{theorem}
Hence, $||T_W||^{-1}$ is the threshold for the final size of the
epidemic on a $G(W)$, i.e.\ if $\beta \leq ||T_W||^{-1}$ there is
a negligible fraction of removed nodes, whereas if $\beta>
||T_W||^{-1}$ there is a non-zero fraction of removed nodes. We
are now going to give explicit computations in the case of
scale-free networks.

Let us examine the case where $$W(x,y)=\frac{W(x)W(y)}{\int_0^1
W(u)du}\:,$$ and such that for $X$ uniformly distributed on
$[0,1]$, the variable $W(X)$ follows a Pareto (power-law)
distribution, i.e.
\begin{equation}\label{eq-power-law}
\P(W(X)\in(t,t+dt))=\frac{\gamma}{(1+t)^{1+\gamma}},\:
t\in[0,\infty)\:.
\end{equation}
This yields $W(x)=(1-x)^{-1/\gamma}-1$ and if we assume that
$\gamma>2$, we have
$$W(x,y)=\frac{\gamma-1}{n}((1-x)^{-1/\gamma}-1)((1-y)^{-1/\gamma}-1)\:.$$
Let us now go back to the operator $T_W$, our goal is to compute
the $L_2$-module of
$$(T_Wf)(x)=\int_0^1 n
\frac{\gamma-1}{n}((1-x)^{-1/\gamma}-1)((1-y)^{-1/\gamma}-1)
f(y)dy\:.$$ It is not difficult to see that the maximum is reached
for $f(y)=\frac{((1-y)^{-1/\gamma}-1)}{\int_0^1 W(u)du}$ and thus
$$
||T_W||=(\gamma-1)\int_0^1W(y)^2dy=\frac{2}{\gamma-2}\:.
$$
Applying the above results for this specific $W-$graph we see that
there will a large outbreak if we have a giant component in the
corresponding graph that is to say $\frac{2\beta}{\gamma-2}>1$.
Hence if $\beta
>\frac{\gamma-2}{2}$ then the epidemic will eventually reach a
proportion $\tau>0$ of the total population. As a by product of
\cite[Theorem 6.2]{BJR05}, it turns out that $\tau$, is related to
the solution of the following functional equation
\begin{equation}\label{eq-functional}
f=1-e^{T_Wf}\:.
\end{equation}
We omit this last computation as it is a bit tedious and does not
bring any further insight into the model. As previously we can
conclude that if $\beta
>\frac{\gamma-2}{2}$ then the final size of the epidemic is
bounded below by $(1+o(1))\tau^2n$ where $\tau$ is related to the
solution of the functional equation (\ref{eq-functional}).

Finally note that if $T$ is a random variable with probability
distribution given by (\ref{eq-power-law}) then the epidemic
threshold $\frac{\gamma-2}{2}$ corresponds to the ratio
$\E(T)/\E(T^2)$. Moreover this gives a stronger result than the
one derived in paragraph \ref{sec:powerlaw} since
$\frac{\gamma-2}{2}<\gamma-1$ where $\frac{1}{\gamma-1}$ is the
average degree $d$.

\section{Conclusion}\label{sec:summary}

Probabilistic methods and tools offer a powerful set of analytical
techniques to understand the spread of epidemics. Such techniques
were used in this paper to gain further insight into models, which
typically have been investigated through mean-field approximations
and simulation studies. Let us recapitulate our key results and
pinpoint some further directions of research.

We derived a threshold for a small outbreak and showed that it is
indeed close to the best possible in general, as it has been be
demonstrated for the example of the star-shaped network. We then
took advantage of results characterising the giant component for
various families of graphs to give a lower bound to the number of
nodes ultimately removed. The Reed-Frost model represents a
starting point and it would be useful to extend our analysis to
exponentially distributed and general infectious periods. Finally,
we intend to pursue this analysis for other epidemic models and
other classes of topologies.

\end{document}